
\documentstyle{amsppt}
\baselineskip18pt
\magnification=\magstep1
\pagewidth{30pc}
\pageheight{45pc}
\hyphenation{co-deter-min-ant co-deter-min-ants pa-ra-met-rised
pre-print pro-pa-gat-ing pro-pa-gate
fel-low-ship Cox-et-er dis-trib-ut-ive}
\def\leaderfill{\leaders\hbox to 1em{\hss.\hss}\hfill}
\def\A{{\Cal A}}

\def\Sy{{\Cal S}}

\def\afn{{\text {\bf a}}}

\def\idest{i.e.,\ }

\def\a{{\alpha}}
\def\be{{\beta}}
\def\g{{\gamma}}

\def\d{{\delta}}

\def\e{{\varepsilon}}

\def\th{{\theta}}

\def\k{{\kappa}}
\def\l{{\lambda}}

\def\s{{\sigma}}
\def\t{{\tau}}

\def\bc{{\bold c}}

\def\BB{{\bold B}}
\def\b0{\text{\bf 0}}

\def\ra{{\ \longrightarrow \ }}

\def\aut{\text{\rm \, Aut}}

\def\real{{\Bbb R}}
\def\complex{{\Bbb C}}
\def\zed{{\Bbb Z}}

\def\Im{\text{\rm Im}}

\def\boxit#1{\vbox{\hrule\hbox{\vrule \kern3pt
\vbox{\kern3pt\hbox{#1}\kern3pt}\kern3pt\vrule}\hrule}}
\def\rabbit{\vbox{\hbox{\kern0pt
\vbox{\kern0pt{\hbox{---}}\kern3.5pt}}}}

\def\tableau#1{
        \hbox {
                \hskip -10pt plus0pt minus0pt
                \raise\baselineskip\hbox{
                \offinterlineskip
                \hbox{#1}}
                \hskip0.25em
        }
}

\def\tabCol#1{
\hbox{\vtop{\hrule
\halign{\strut\vrule\hskip0.5em##\hskip0.5em\hfill\vrule\cr\lower0pt
\hbox\bgroup$#1$\egroup \cr}
\hrule
} } \hskip -10.5pt plus0pt minus0pt}

\def\CR{
        $\egroup\cr
        \noalign{\hrule}
        \lower0pt\hbox\bgroup$
}



\def\blank#1#2{
\hbox to #1{\hfill \vbox to #2{\vfill}}
}


\def\strut{\vrule height10pt depth5pt width0pt}

\topmatter
\title
Categories arising from tabular algebras
\endtitle

\author R.M. Green \endauthor
\affil 
Department of Mathematics and Statistics\\ Lancaster University\\
Lancaster LA1 4YF\\ England\\
{\it  E-mail:} r.m.green\@lancaster.ac.uk\\
\endaffil

\abstract
We continue the investigation of tabular algebras with trace 
(a certain class of associative $\zed[v, v^{-1}]$-algebras equipped with
distinguished bases) by determining the extent to which the tabular
structure may be recovered from a knowledge of the structure
constants.  This problem is equivalent to understanding a certain
category (the category of table data associated to a tabular algebra) 
which we introduce.  The main result is that this category is
equivalent to another category (the category of based posets
associated to a tabular algebra) whose structure we describe
explicitly.
\endabstract


\thanks
\noindent 2000 {\it Mathematics Subject Classification.} 16B50
\endthanks

\endtopmatter



\centerline{\bf To appear in the Glasgow Mathematical Journal}

\head Introduction \endhead

Tabular algebras with trace were introduced by the author in 
\cite{{\bf 4}} as a
class of algebras over the ring $\A := \zed[v, v^{-1}]$.  They are by
definition equipped with a tabular basis that is described in terms
of a ``table datum'' and is required to satisfy various axioms.  There
are many natural examples of tabular bases, including the
Kazhdan--Lusztig bases for certain Hecke algebras \cite{{\bf 8}}, 
the diagram bases of the Brauer algebra or Jones' annular algebra 
\cite{{\bf 7}}, and the IC bases of various kinds of Temperley--Lieb 
algebra \cite{{\bf 5}}.  Tabular algebras also provide a convenient starting
point from which to study cellular algebras in the sense of Graham and
Lehrer \cite{{\bf 3}}.  Cellular algebras are of considerable
interest in representation theory, and several constructions of cell data
for specific algebras in the literature may be unified by using tabular 
algebras \cite{{\bf 4}, Theorem 2.1.1}.  Tabular algebras with trace are
intriguing objects in their own right because all the natural examples
(in fact, all the examples currently known to the author) have additional
properties not required of them by the defining axioms, notably positivity of
structure constants.

It is natural to wonder to what extent the table datum is determined by the
tabular basis, particularly as the definition of a tabular algebra
looks superficially complicated.  More precisely, given a tabular
algebra $A$ with distinguished basis $\BB$ and trace $\t$ satisfying axioms
(A1)--(A5) (see \S2.1), we wish to classify all possible table data
compatible with the basis.  A convenient way to do this is by using the
technique of categorification, thus replacing the sets $A$ and $\BB$ by a 
suitable category, ${\Cal D}(A, \BB)$.  The objects and morphisms of
${\Cal D}(A, \BB)$ are defined in terms of the possible table data for 
$(A, \BB)$.  Understanding the possible table data for $(A, \BB)$ is then
equivalent to understanding the structure of the category
${\Cal D}(A, \BB)$.

Our main result (Theorem 3.1.6) is that the category ${\Cal D}(A,
\BB)$ is equivalent to another category ${\Cal P}(A, \BB)$, whose 
objects and morphisms can be easily and explicitly described in terms 
of ``based posets'', which we introduce.  
This solves the problem of understanding the structure of
${\Cal D}(A, \BB)$.  We also show that the algebra automorphisms
of $A$ that fix $\BB$ setwise may be understood in
terms of these categories, and we show how these may be
computed in typical cases.

Before we can define based posets and state the main results, it is 
necessary to develop some
elementary theory regarding matrix rings over table algebras,
and their automorphisms.  This is the subject of \S1.  We recall the
definition of tabular algebras from \cite{{\bf 4}} in \S2.
In \S3, we introduce and study based posets.  In \S4, we illustrate
some of the ideas of this paper using the Brauer algebra as a worked
example; the reader unfamiliar with tabular algebras may prefer to
look at \S4.2 before reading \S2.

The results of this paper are interesting largely because of their
applications to representation theory.  For example, it is possible to 
define combinatorially 
a class of ``standard modules'' for a tabular algebra $(A, \BB)$ 
in terms of the table datum, and the results of this paper can be used to 
show that the
class of modules so obtained depends only on the pair $(A, \BB)$, and not on
the table datum chosen.  Applications such as these will be explored
in a sequel to this paper, where it will be shown that the extended affine
Hecke algebra of type $\widetilde{A}_n$ equipped with its
Kazhdan--Lusztig basis is a tabular algebra with trace, and that the 
standard modules agree with the geometrically defined standard modules 
appearing in the work of Lusztig \cite{{\bf 9}}.

\head 1. Based rings and their automorphism groups \endhead

\subhead 1.1 Table algebras \endsubhead

We begin by recalling the definition of a table algebra, which is a
generalization of the integral group ring of an arbitrary group.

\definition{Definition 1.1.1}
A table algebra is a pair $(\Gamma, B)$, where $\Gamma$ is an associative 
unital $R$-algebra for some $\zed \leq R \leq \complex$ 
and $B = \{b_i : i \in I\}$ is a distinguished basis for
$\Gamma$ such that $1 \in B$, satisfying the following three axioms:

\item{(T1)}{The structure constants of $\Gamma$ with respect to the basis
$B$ lie in $\real^+$, the nonnegative real numbers.}
\item{(T2)}{There is an algebra anti-automorphism $\bar{\ }$ of $\Gamma$ whose
square is the identity and that has the property 
that $b_i \in B \Rightarrow \overline{b_i} \in
B$.  (We define $\overline{i}$ by the condition $\overline{b_i} =
b_{\bar{i}}$.)}
\item{(T3)}{Let $\k(b_i, a)$ be the coefficient of $b_i$ in $a \in \Gamma$.
Then there is a function $g: B \times B \ra \real^+$ satisfying $$
\k(b_m, b_i b_j) = g(b_i, b_m) \k(b_i, b_m \overline{b_j})
,$$ where $g(b_i, b_m)$ is independent of $j$, for all $i, j, m$.
}
\enddefinition

\remark{Remark 1.1.2}
Table algebras first appeared in the work of Arad and Blau \cite{{\bf 1}} 
in the case where $\Gamma$ is commutative and $B$ is finite.  All table
algebras in this paper will be {\it normalized}, meaning that the
structure constants of $\Gamma$ with respect to the basis $B$ will be
(nonnegative) integers and the function $g$ in axiom (T3) sends all
pairs of basis elements to $1$.  These conditions are reminiscent of
Sunder's discrete hypergroups \cite{{\bf 11}}, and they are clearly
satisfied when $B = G$ is any group and $\Gamma = \zed G$.
\endremark

\proclaim{Lemma 1.1.3}
Let $(\Gamma, B)$ be a normalized table algebra. 
\item{\rm (i)}{The linear map $$t : a \ra
\kappa(1, a)$$ is a trace on $\Gamma$ 
(that is, $t(xy) = t(yx)$ for all $x, y \in \Gamma$).}
\item{\rm (ii)}{Let $a = \sum_{i \in I} z_i b_i$ where $b_i \in B$ and $z_i
\in \zed$.  Then $t(a\bar{a}) = 1$ if and only if $a = \pm b_i$ for
some $i$.}
\endproclaim

\demo{Proof}
Part (i) is an easy consequence of axiom (T3), which shows that
$\kappa(1, b_i b_j)$ = $\d_{i \bar{j}}$.  Using this observation, we
see that $t(a \bar{a}) = \sum_{i \in I} z_i^2$, which proves (ii).
\qed\enddemo

\definition{Definition 1.1.4}
Let $(\Gamma, B)$ be a normalized table algebra.  A basis element $b
\in B$ is said to be {\it grouplike} if $b \bar{b} = \bar{b} b = 1$.
\enddefinition

The notion of a grouplike element is similar to Arad and Blau's notion
of an ``irreducible'' element, but this needs to be restated in our
context so that we can deal with the case where $B$ is infinite.


\subhead 1.2 Based rings \endsubhead

\definition{Definition 1.2.1}
A {\it based ring} is a pair $(A, B)$, where $A$ is a unital $\zed$-algebra
with free $\zed$-basis $B$ and nonnegative structure constants.  A
homomorphism $\phi : (A, B) \ra (A', B')$ of based rings is a
homomorphism of abstract $\zed$-algebras $\phi : A \ra A'$ such that $\phi(b)
\in B' \cup \{0\}$ for all $b \in B$.  Isomorphisms,
automorphisms etc. of based rings are defined analogously.  
\enddefinition

Clearly normalized table algebras are examples of based rings.

\proclaim{Lemma 1.2.2}
Let $(\Gamma, B)$ be a normalized table algebra and let $b, g \in
B$ with $b$ grouplike.  
\item{\rm (i)}{We have $bg \in B$ and $gb \in B$.}
\item{\rm (ii)}{The $\zed$-linear map sending $g' \mapsto \bar{b} g'
b$ for all $g' \in B$ is a based ring automorphism of $(\Gamma, B)$.}
\endproclaim

\demo{Proof}
Claim (ii) is immediate from (i), since $\bar{b}$ is grouplike if and
only if $b$ is, so it remains to prove (i).

Let $t$ be the trace of Lemma 1.1.3.  We observe that $$
t((gb)(\overline{gb})) = t(gb \bar{b} \bar{g}) = t(g\bar{g}) = 1
,$$ where the second equality uses the fact that $b$ is grouplike, and
the third equality uses Lemma 1.1.3 (ii).  Lemma 1.1.3 (ii) and the assumption
that $(\Gamma, B)$ is normalized show that $gb
\in B$.  To prove the other half of (i), note that $t((bg)(\overline{bg}))
= t((\overline{bg})(bg))$ by Lemma 1.1.3 (i) and then proceed as before.
\qed\enddemo

The main example of a based ring that is of interest for our purposes is the
following.

\definition{Definition 1.2.3}
Let $(\Gamma, B)$ be a normalized table algebra.  
The based ring $M_{n, \Gamma, B}$ is the ring of $n \times n$
matrices over the ring $\Gamma$, equipped with distinguished basis
consisting of all elements $e_{ij} \otimes b$, where $e_{ij}$ is a
matrix unit and $b \in B$.  

We call the elements $e_{ii} \otimes 1$ (for $1 \leq i \leq n$) 
{\it distinguished idempotents}.
\enddefinition

It is trivial to check that $M_{n, \Gamma, B}$ is indeed a based
ring.

\proclaim{Lemma 1.2.4}
Let $(\Gamma, B)$ be a normalized table algebra, and let $\bar{\ }$
be the table algebra anti-automorphism.  Let $\a$ be an automorphism
of $M_{n, \Gamma, B}$ (as a based ring).
\item{\rm (i)}{The map $\a$ sends distinguished idempotents of $M_{n,
\Gamma, B}$ to distinguished idempotents.}
\item{\rm (ii)}
{The $\zed$-linear map $*$ that sends $e_{ij} \otimes b$ to $(e_{ij}
\otimes b)^* := e_{ji} \otimes \bar{b}$ is an anti-automorphism 
of $M_{n, \Gamma, B}$ that commutes with $\a$.}
\item{\rm (iii)}
{If $M_{n', \Gamma', B'}$ is a based ring with anti-automorphism 
$*'$ (as in (ii)) and $\phi : M_{n, \Gamma, B} \ra M_{n', \Gamma', B'}$ 
is an isomorphism of based rings, then $*' = \phi \circ * \circ \phi^{-1}$.}
\endproclaim

\demo{Proof}
To prove (i), we note that $\a$ preserves the identity element, which is
expressed in terms of the based ring basis as $$
1 = \sum_{i = 1}^n (e_{ii} \otimes 1)
.$$  As an automorphism of based rings, $\a$ permutes the basis
elements, and (i) follows.

It is easy to check that the map $*$ is an anti-automorphism of based
rings.  

Consider two basis elements 
$e_{ij} \otimes b$ and $e_{kl}
\otimes b'$ of the based ring.  The only way a distinguished
idempotent can occur with nonzero coefficient in the product
$(e_{ij} \otimes b)(e_{kl} \otimes b')$
is if $j = k$, $i = l$ and $1$ occurs in the product $bb'$.
The last of these conditions happens if and only if $b' = \bar{b}$ by
axiom (T3).  If all these conditions hold, we have $e_{kl} \otimes b' =
(e_{ij} \otimes b)^*$, and the only distinguished idempotent occurring
in the product is $e_{ii} \otimes 1$, which occurs with coefficient
$1$ since $(\Gamma, B)$ is normalized.  This characterizes $*$ in terms of the
structure constants and distinguished idempotents, and (iii) follows.

Now consider the basis element $\a(e_{ij} \otimes b)$.  Since $\a$
permutes the distinguished idempotents by (i), we may apply $\a$ to the
equation $$
(e_{ij} \otimes b)(e_{ji} \otimes \bar{b}) = \sum c_{p, q, b'} (e_{pq}
\otimes b')
$$ and argue as in the previous paragraph to show that $$
\a(e_{ji} \otimes \bar{b}) = \a((e_{ij} \otimes b)^*)
= (\a(e_{ij} \otimes b))^*
.$$  Claim (ii) follows by linearity.
\qed\enddemo

\subhead 1.3 Automorphisms of $M_{n, \Gamma, B}$ \endsubhead

In \S1.3, we take a closer look at the based rings $M_{n, \Gamma, B}$ of 
Definition 1.2.3.
The following lemma shows how the based ring isomorphism type of such a ring
is controlled by the data $n$, $\Gamma$ and $B$.

\proclaim{Lemma 1.3.1}
If $\a : M_{n, \Gamma, B} \ra M_{n', \Gamma', B'}$ is an isomorphism
of based rings, then $n = n'$ and $(\Gamma, B) \cong (\Gamma',
B')$ as based rings.
\endproclaim

\demo{Proof}
Since $\a$ is an isomorphism of based rings, it induces a bijection
between the bases of each based ring.  Arguing as in the proof of
Lemma 1.2.4 (i), we see that $\a$ sends distinguished idempotents to
distinguished idempotents, so that in particular we have $n = n'$.  

Let $e$ be any distinguished idempotent in $M_{n, \Gamma, B}$ and
let $\phi_e$ be the $\zed$-linear map from $(\Gamma, B)$ to $M_{n, \Gamma,
B}$ for which $\phi_e(b) = e \otimes b$ for all $b \in B$.  It is
clear that $\phi_e$ is a monomorphism of based rings whose image is $e
M_{n, \Gamma, B} e$, which shows that the isomorphism type of
$(\Gamma, B)$ as a based ring is determined by that of $M_{n,
\Gamma, B}$.  The conclusion follows.
\qed\enddemo

The classification of automorphisms of $M_{n, \Gamma, B}$ is more
interesting than the proof of Lemma 1.3.1 suggests.  This is
due to the presence of what we call ``twisted'' isomorphisms (defined below) 
which may not 
send elements of the form $e_{ij} \otimes 1$ to elements of the form
$e_{kl} \otimes 1$ if $i \ne j$.

\definition{Definition 1.3.2}
Let $\a : M_{n, \Gamma, B} \ra M_{n', \Gamma', B'}$ be a homomorphism
of based rings.  If there exist a map $\s : \{1, \ldots, n\} \ra 
\{1, \ldots, n'\}$ and a homomorphism of based rings 
$\psi : (\Gamma, B) \ra (\Gamma', B')$ 
such that $\a(e_{ij} \otimes b) = e_{\s(i) \s(j)} \otimes 
\psi(b)$ for all $1 \leq i, j \leq n$ and $b \in B$ then we call $\a$ an
{\it untwisted} homomorphism.  Otherwise, we call $\a$ a {\it twisted}
homomorphism.
\enddefinition

Lemma 1.3.1 has the following

\proclaim{Corollary 1.3.3}
If $M_{n, \Gamma, B}$ and $M_{n', \Gamma', B'}$ are isomorphic as based rings,
then they are isomorphic by an untwisted isomorphism. \qed
\endproclaim

\definition{Definition 1.3.4}
Consider the based ring $M_{n, \Gamma, B}$.

We associate to the sequence $(b_1, b_2, \ldots, b_n)$ of grouplike
elements of $B$ the based ring automorphism $\be = \be(b_1, b_2, \ldots,
b_n)$ of $M_{n, \Gamma, B}$.  This is defined to send the element
$X \in M_{n, \Gamma, B}$ to $G^{-1} X G$, where $$
G = \sum_{i = 1}^n e_{ii} \otimes b_i
.$$  (This makes sense by Lemma 1.2.2.)  The automorphism $\be$ will be
twisted unless all the $b_i$ are equal.

If $w$ is a permutation in the symmetric group $\Sy(n)$, we define an
untwisted based ring automorphism, $\s_w$, of $M_{n, \Gamma, B}$ that sends
$e_{ij} \otimes b$ to $e_{w^{-1}i, w^{-1}j} \otimes b$.

If $\e$ is a based ring automorphism of $(\Gamma, B)$ (not
necessarily one of the form given in Lemma 1.2.2 (ii)), we define an
untwisted based ring automorphism $\psi_\e$ of $M_{n, \Gamma, B}$ by
$\psi_\e (e_{ij} \otimes b) = e_{ij} \otimes \e(b)$.
\enddefinition

\proclaim{Proposition 1.3.5}
The group $\aut_B(M_{n, \Gamma, B})$ of based ring automorphisms 
of $M_{n, \Gamma, B}$ is generated by the automorphisms of the form 
$\be(b_1, b_2, \ldots, b_n)$, $\s_w$ and $\psi_\e$ as 
given in Definition 1.3.4.
\endproclaim

\demo{Proof}
Let $\a$ be such an automorphism.  By Lemma 1.2.4 (i), $\a$ permutes
the distinguished idempotents, so by applying a suitable automorphism
$\s_w$, we may assume without loss of generality that $\a(e_{ii}
\otimes 1) = e_{ii} \otimes 1$ for all $1 \leq i \leq n$.  

With the above assumption, it follows 
that for any $1 \leq k, l \leq n$ and $b \in B$,  
$\a(e_{kl} \otimes b)$ is of the form $e_{kl} \otimes b'$.  This is
because $$
(e_{kk} \otimes 1)(e_{kl} \otimes g)(e_{ll} \otimes 1) = e_{kl}
\otimes g
$$ for any $g \in B$.

For $2 \leq i \leq n$, we define the element $b'_i \in B$ by the
condition $\a(e_{i, i-1} \otimes 1) = e_{i, i-1} \otimes b'_i$.

We claim that $b'_i$ is grouplike.  Recalling the map $*$ from Lemma 1.2.4
(ii), we note that $(e_{i, i-1} \otimes 1)(e_{i, i-1} \otimes 1)^* = e_{ii}
\otimes 1$.  Since $\a$ commutes with $*$ and fixes the distinguished
idempotents, we must have $(e_{i, i-1} \otimes b'_i)(e_{i-1, i} \otimes
\bar{b'_i}) = e_{ii} \otimes 1$, which implies that $b'_i \bar{b'_i} = 1$.  A
similar argument establishes that $\bar{b'_i} b_i' = 1$, so that $b_i'$ is
grouplike.

Let $b_1 = 1$ and define, for $2 \leq i \leq n$, $b_i := b'_i b'_{i-1}
\cdots b'_2$.  Let $\be$ be the automorphism $\be(b_1, b_2, \ldots,
b_n)$ of Definition 1.3.4.  (Note that $\be$ commutes with $*$; this
follows easily from the definition since the matrix $G$ in Definition
1.3.4 is diagonal.)  A routine matrix calculation shows that
$\be(\a(e_{i, i-1} \otimes 1)) = e_{i, i-1} \otimes 1$ for all $2 \leq
i \leq n$.  Since $\be$ and $\a$ commute with $*$, the map $\be \a$ also
fixes elements of the form $e_{i-1, i} \otimes 1$, and since $\be \a$
is an algebra homomorphism, it fixes all elements of the form $e_{ij}
\otimes 1$ for $1 \leq i, j \leq n$.
We may therefore assume for the rest of the proof that $\a$
fixes all elements $e_{ij} \otimes 1$.

We now see that $\a$ is determined by its values on $e_{11} \otimes
b$, because $$e_{ij} \otimes b = (e_{i1} \otimes 1)(e_{11} \otimes
b)(e_{1j} \otimes 1)
.$$  By the second paragraph of the proof, $\a$ must send $e_{11} \otimes
b$ to $e_{11} \otimes b'$ for some $b'$.  As noted in the proof of
Lemma 1.3.1, the $\zed$-linear map sending $b$ to $e_{11} \otimes b$
is a based ring monomorphism, so it follows that there is a
based ring automorphism $\e$ of $(\Gamma, B)$ such that $\a(e_{11}
\otimes b) = e_{11} \otimes \e(b)$.  In other words, $\a$ is equal to
$\psi_\e$, completing the proof.
\qed\enddemo

\head 2. Tabular algebras and their based rings \endhead

We now recall from \cite{{\bf 4}} 
the definition of a tabular algebra and its associated
table datum.  The goal of this paper is to understand the extent to
which the table datum is determined by the tabular basis.

\subhead 2.1 Tabular algebras \endsubhead

We start by recalling the definition of the $\afn$-function, which is
due to Lusztig.

\definition{Definition 2.1.1}
Let $\A$ be the ring of Laurent polynomials $\zed[v, v^{-1}]$, let
$A$ be an $\A$-algebra and let $\BB$ be an $\A$-basis of $A$.
For $X, Y, Z \in \BB$, we define the structure constants $g_{X, Y, Z} \in
\A$ by the formula $$
X Y = \sum_Z g_{X, Y, Z} Z
.$$  The $\afn$-function is defined by $$
\afn(Z) = 
\max_{X, Y \in \BB} \deg(g_{X, Y, Z})
,$$ where the degree of a Laurent polynomial is taken to be 
the highest power of $v$ occurring with nonzero coefficient.  We
define $\g_{X, Y, Z} \in \zed$ to be the coefficient of $v^{\afn(Z)}$ in $g_{X,
Y, Z}$; this will be zero if the bound is not achieved.
\enddefinition

\definition{Definition 2.1.2}
A {\it tabular algebra} is an
$\A$-algebra $A$, together with a table datum 
$$(\Lambda, \Gamma, B, M, C, *)$$ satisfying axioms (A1)--(A3) below.

\item{(A1)}
{$\Lambda$ is a finite poset.  For each $\l \in \Lambda$, 
$(\Gamma(\l), B(\l))$ is
a normalized table algebra over $\zed$ and
$M(\l)$ is a finite set.  The map $$
C : \coprod_{\l \in \Lambda} \left( M(\l) \times B(\l) \times M(\l)
\right) \rightarrow A
$$ is injective with image an $\A$-basis of $A$.  We assume
that $\Im(C)$ contains a set of mutually orthogonal idempotents 
$\{1_\e : \e \in {\Cal E}\}$ such that 
$A = \sum_{\e, \e' \in {\Cal E}} (1_\e A 1_{\e'})$ and such that for each
$X \in \Im(C)$, we have $X = 1_\e X 1_{\e'}$ for some $\e, \e' \in
{\Cal E}$.
A basis arising in
this way is called a {\it tabular basis}.  
}
\item{(A2)}
{If $\l \in \Lambda$, $S, T \in M(\l)$ and $b \in B(\l)$, we write
$C(S, b, T) = C_{S, T}^{b} \in A$.  
Then $*$ is an $\A$-linear involutory anti-automorphism 
of $A$ such that
$(C_{S, T}^{b})^* = C_{T, S}^{\overline{b}}$, where $\bar{\ }$ is the
table algebra anti-automorphism of $(\Gamma(\l), B(\l))$.
If $g \in \complex(v) \otimes_\zed \Gamma(\l)$ is such that 
$g = \sum_{b_i \in B(\l)} c_i b_i$ for some scalars $c_i$ 
(possibly involving $v$), we write
$C_{S, T}^g \in \complex(v)\otimes_\A A$ 
as shorthand for $\sum_{b_i \in B(\l)} c_i C_{S, T}^{b_i}$.  We write
$\bc_\l$ for the image under $C$ of $M(\l) \times B(\l) \times M(\l)$.}
\item{(A3)}
{If $\l \in \Lambda$, $g \in \Gamma(\l)$ and $S, T \in M(\l)$ then for all 
$a \in A$ we have $$
a . C_{S, T}^{g} \equiv \sum_{S' \in M(\l)} C_{S', T}^{r_a(S', S) g}
\mod A(< \l),
$$ where  $r_a (S', S) \in \Gamma(\l)[v, v^{-1}] = \A \otimes_\zed
\Gamma(\l)$ is independent of $T$ and of $g$ and $A(< \l)$ is the
$\A$-submodule of $A$ generated by the set $\bigcup_{\mu < \l} \bc_\mu$.}

A {\it tabular algebra with trace} is a tabular algebra that also
satisfies conditions (A4) and (A5) below.

\item{(A4)}{Let $K = C_{S, T}^b$, $K' = C_{U, V}^{b'}$ and 
$K'' = C_{X, Y}^{b''}$ lie in $\Im(C)$.  Then the
maximum bound for $\deg(g_{K, K', K''})$
in Definition 2.1.1 is achieved if and only if $X = S$, $T = U$, $Y =
V$ and $\kappa(b'', bb') \ne 0$ (where $\kappa$ is as in axiom (T3)).  
If these conditions all hold and
furthermore $b = b' = b'' = 1$, we require $\g_{K, K', K''} = 1$.}
\item{(A5)}{There exists an $\A$-linear function $\t : A \ra \A$
(the {\it tabular trace}), such that $\t(x) = \t(x^*)$ for all $x \in
A$ and $\t(xy) = \t(yx)$ for all $x, y \in A$, that has the 
property that for every
$\l \in \Lambda$, $S, T \in M(\l)$, $b \in B(\l)$ and $X = C_{S,
T}^b$, we have $$
\t(v^{\afn(X)} X) = 
\cases 1 \mod v^{-1} \A^- & \text{ if } S = T \text{ and } b = 1,\cr
0 \mod v^{-1} \A^- & \text{ otherwise.} \cr
\endcases
$$  Here, $\A^- := \zed[v^{-1}]$.  We call the elements $C_{S, S}^1$
{\it distinguished involutions}.}
\enddefinition

\remark{Remark 2.1.3}
The idempotent condition in axiom (A1) ensures that $\afn(Z)$ is
always defined.
\endremark

Tabular algebras are so called because they are an amalgamation of
table algebras and cellular algebras in the sense of Graham and Lehrer
\cite{{\bf 3}}.  Axioms (A1)--(A3) are modelled on the axioms for a cellular
algebra.  

In this paper, we will only be concerned with tabular
algebras with trace; this class of examples includes all the examples
mentioned in the introduction.  Our goal is to show that although the
axioms in Definition 2.1.2 seem complicated, one can recover the table
datum ``up to isomorphism'' (in a sense that will be made precise)
from the structure constants of the tabular basis.  
Another way to state our aim is by the following question.

\proclaim{Question 2.1.4}
Given a tabular algebra $A$ with trace and with tabular basis $\BB$, 
to what extent can we recover the map $C$?
\endproclaim

Question 2.1.4 can be viewed as a question about categories, as we
will explain in \S2.4.

It is too much to hope to recover the table datum from the structure of
$A$ as an abstract algebra, as can be seen from the following result.

\proclaim{Theorem 2.1.5 (Hertweck \cite{{\bf 6}})}
There exist finite groups $G$ and $H$ with $G \not\cong H$ and $\zed G
\cong \zed H$ (so that $\A G = \A \otimes_\zed G \cong \A \otimes_\zed
H = \A H$).
\endproclaim

The relevance to tabular algebras is as follows.  Since $G$ is a group,
$\A G$ is a tabular algebra with trace: take $\Lambda$ to consist of a
single element $\l$, $\Gamma(\l) = \zed G$, $B(\l) = G$, 
$M(\l)$ to be a single element $m$, $C(m, b, m) = b$ and $*$ to be the
linear extension of inversion.  We can take $\t(x)$ to be the
coefficient of the identity element in $x$.  Hertweck's theorem then
shows that the isomorphism type of $(\Gamma(\l), B(\l))$ as a based
ring cannot be recovered from the isomorphism type of $\A G$ as an
abstract algebra.

\subhead 2.2 Based rings arising from tabular algebras \endsubhead

In \cite{{\bf 4}, \S3}, asymptotic versions of tabular algebras with trace
are constructed, using methods from \cite{{\bf 10}}.  These asymptotic
algebras are based rings in the sense of \S1.2.  They will be useful
in answering Question 2.1.4 since it will turn out that we can recover
information about the tabular algebra by studying the associated based
ring.

\definition{Definition 2.2.1}
Let $A$ be a tabular algebra with trace, and maintain the usual notation.
Define $\widehat{X} := v^{-\afn(X)} X$ for any
tabular basis element $X \in \Im(C)$.  
The free $\A^-$-module $A_\l^-$ is
defined to be generated by the elements $\{\widehat{X} : X \in
\bc_\l\}$.  We set $t_X$ to be the image of $\widehat{X}$ in $$
A_\l^\infty := {{A_\l^-} \over {v^{-1} A_\l^-}}
.$$  The latter is a $\zed$-algebra with basis $\{t_X : X \in
\bc_\l\}$ and structure constants $$
t_X t_{X'} = \sum_{X'' \in \bc} \g_{X, X', X''} t_{X''}
,$$ where the $\g_{X, X', X''} \in \zed$ are as in Definition 2.1.1.
We also set $$
A^\infty := \bigoplus_{\l \in \Lambda} A_\l^\infty
;$$ this is a $\zed$-algebra with basis $\{t_X : X \in \Im(C)\}$.
\enddefinition

We will call the ring $A^\infty$ the {\it based ring associated to the
tabular algebra} $A$.  This terminology is justified by the following lemma.

\proclaim{Lemma 2.2.2}
Let $A$ be a tabular algebra with trace, and maintain the usual
notation.  
\item{\rm (i)}{The algebra $A_\l^\infty$ with basis $\{t_X : X \in \bc_\l\}$ is
isomorphic as a based ring to $M_{|M(\l)|, \Gamma(\l), B(\l)}$.  The
isomorphism may be chosen to
identify $t_X$ with $e_{ST} \otimes b$, where $X = C(S, b, T)$ and
$M(\l)$ is identified with the set $\{1, 2, \ldots, |M(\l)|\}$.}
\item{\rm (ii)}{The algebra $A^\infty$ with basis $\{t_X : X \in
\Im(C)\}$ is a based ring.}
\endproclaim

\demo{Proof}
Part (i) follows from \cite{{\bf 4}, Theorem 3.2.4 (i)} and its proof.
Part (ii) is immediate from part (i) and the definition of $A^\infty$.
\qed\enddemo

\subhead 2.3 Reduced tabular algebras \endsubhead

It is clear from axiom (A3) of a tabular algebra that if $A$ is a
tabular algebra with table datum $(\Lambda, \Gamma, B, M, C, *)$ then
we may refine the partial order $\leq$ on $\Lambda$ to a larger
partial order without disturbing any of the axioms.  However, this
extra freedom turns out to be inconvenient for our purposes in this
paper since it obfuscates some of the symmetry properties of the cell
datum.  For this reason, we introduce the notion of a reduced tabular
algebra, for which the partial order on $\Lambda$ is as small as possible.
Most of our results concern reduced tabular algebras, but there is no
loss of generality in assuming that a tabular algebra is reduced.

\definition{Definition 2.3.1}
Let $A$ be a tabular algebra with table datum $(\Lambda, \Gamma, B, M,
C, *)$, where $\Lambda$ is ordered by $\leq$.  Let $\Lambda'$ be a
poset with the same underlying set as $\Lambda$, partially ordered by
$\leq'$, and write $\Lambda' \preceq_P \Lambda$ if $\leq$ is a
refinement of $\leq'$.  If $(\Lambda', \Gamma, B, M, C, *)$ is a table
datum for $A$ for some $\Lambda' \prec_P \Lambda$ then $A$ and its
table datum are said to be non-reduced; otherwise, $A$ and its table
datum are said to be {\it reduced}.
\enddefinition

One of the advantages of reduced tabular algebras with trace is that the poset
$\Lambda$ may be recovered up to isomorphism from the tabular basis.

\definition{Definition 2.3.2}
Let $A$ be a tabular algebra.  
If $X$ and $X'$ are tabular basis elements, we say that 
$X' \preceq X$ if $X'$ appears with nonzero coefficient in $KXK'$ for 
some tabular basis elements $K, K'$.  The relation $\preceq$ on
the tabular basis $\Im(C)$ is defined to be the transitive extension
of this relation; it is reflexive by axiom (A1).
\enddefinition

The following result is the first step towards recovering the table
datum of a tabular algebra with trace from the structure constants.

\proclaim{Proposition 2.3.3}
Let $A$ be a reduced tabular algebra with trace and table datum 
$$(\Lambda, \Gamma, B, M, C, *).$$  Let $\preceq$ be as in Definition
2.3.2.  Let $X = C_{S, T}^b \in \bc_\l$ and $X' = C_{U, V}^{b'} \in
\bc_\mu$ be tabular basis elements (where $\l, \mu \in \Lambda$ and 
$\bc_\l$ is as defined in axiom (A2)).  Then $X \preceq X'$ if and
only if $\l \leq \mu$, with
$\l = \mu$ if and only if $X \preceq X' \preceq X$.  It follows that
the tabular basis determines the sets $\bc_\l$ and the isomorphism
type of the poset $\Lambda$.
\endproclaim

\demo{Proof}
It is clear from axiom (A3), its mirror image under $*$ (see 
\cite{{\bf 4}, Remark 1.3.2}) and the definition of $\preceq$ that 
$X \preceq X'$ implies $\l \leq \mu$.  

By \cite{{\bf 4}, Proposition 3.1.3}, we 
find that $\l = \mu$ (and thus $\bc_\l = \bc_\mu$) if and only if 
$X \preceq X' \preceq X$.  To complete the proof of 
the first assertion, it remains to show
that if $\l < \mu$ then $X \preceq X'$.

Assume $\l < \mu$.
Since $A$ is reduced, the partial order $\leq$ is the smallest possible 
partial order compatible with axiom (A3).  The fact that $\Lambda$ is finite
means that there is a chain $$
\l = \l_1 < \l_2 < \cdots < \l_r = \mu
$$ where, for each $1 \leq i < r$, there exist basis elements 
$Y_i \in \bc_{\l_i}$, $Y_{i+1} \in \bc_{\l_{i+1}}$ and $K \in \BB$ such that 
$Y_i$ occurs with nonzero coefficient in the expansion of $KY_{i+1}$.  The
idempotent condition of axiom (A1) shows that $Y_i \preceq Y_{i+1}$, and
we have $Y_1 \preceq Y_r$ by transitivity.  The previous paragraph shows
that $X \preceq Y_1 \preceq X$ and $X' \preceq Y_r \preceq X'$, so that
$X \preceq X'$ as required.

The second assertion now follows from the observation that 
the definition of $\preceq$ depends only on the tabular basis and 
not on any details of the table datum.
\qed\enddemo

\proclaim{Corollary 2.3.4}
Let $A$ be a reduced tabular algebra with trace and let 
$$(\Lambda, \Gamma, B, M, C, *), \quad (\Lambda', \Gamma', B', M', C', *')$$ 
be two table data for $A$ associated to the same tabular basis 
$\BB = \Im C = \Im C'$.  Let $X \in \BB$, and define $\l \in \Lambda$ and 
$\l' \in \Lambda'$ by the conditions $X \in \bc_\l$ and $X \in \bc_{\l'}$.
There is an isomorphism of based rings $$
M_{|M(\l)|, \Gamma(\l), B(\l)} \cong
M_{|M(\l')|, \Gamma(\l'), B(\l')}
.$$
\endproclaim

\demo{Proof}
By Proposition 2.3.3, the set $\bc_\l$ containing $X$ may be
reconstructed from $X$ and the tabular basis.  This enables us to
recover $A_\l^\infty$ from Definition 2.2.1, and the conclusion
follows from Lemma 2.2.2 (i).
\qed\enddemo

\subhead 2.4 Categories arising from table data \endsubhead

Question 2.1.4 can be restated in terms of a certain category that we
now introduce.

\definition{Definition 2.4.1}
Let $A$ be a tabular algebra with trace $\t$ and tabular basis $\BB$.  The
category ${\Cal D}(A, \BB) = {\Cal D}(A, \BB, \t)$ is defined as follows.

\noindent Objects: All elements $(\Lambda, \Gamma, B, M, *)$ for which
there exists $C$ such that $$(\Lambda, \Gamma, B, M, C, *)$$ is a reduced table
datum for $A$ with $\Im(C) = \BB$.

\noindent Morphisms:  Let $(\Lambda, \Gamma, B, M, *)$ and
$(\Lambda', \Gamma', B', M', *')$ be objects of ${\Cal D}(A, \BB)$,
and fix $C$ such that $(\Lambda, \Gamma, B, M, C, *)$ is a reduced table
datum for $A$ with $\Im(C) = \BB$.
The set of morphisms between $(\Lambda, \Gamma, B, M, *)$
and $(\Lambda', \Gamma', B', M', *')$ are the maps $$
(C')^{-1} \circ C : 
\coprod_{\l \in \Lambda} \left( M(\l) \times B(\l) \times M(\l)
\right) \ra 
\coprod_{\l' \in \Lambda'} \left( M'(\l') \times B'(\l') \times M'(\l')
\right)
,$$ where $C'$ is such that $(\Lambda', \Gamma', B', M', C', *')$ is a
reduced table datum for $A$ with $\Im(C') = \BB$, and
composition is given by composition of maps.
\enddefinition

\remark{Remark 2.4.2} It is not clear at this stage that ${\Cal D}(A,
\BB)$ is a well-defined category, because it is not {\it a priori}
obvious that the composition of two morphisms is another morphism or
that the set of morphisms between two objects is independent of the
choice of $C$.  We will resolve this issue in Theorem 3.4.1.
\endremark

Question 2.1.4 is asking for a classification of the morphisms from a given
object in ${\Cal D}(A, \BB)$.  We will achieve this by exhibiting an
equivalence of categories between ${\Cal D}(A, \BB)$ and a category
for which this question is easy to answer.


\head 3. Based posets and their automorphisms \endhead

We now introduce the notion of a based poset, which allows us to
state our main result, Theorem 3.1.6.  Throughout \S3, $A$ will be a
reduced tabular algebra with trace $\t$ and tabular basis $\BB$.

\subhead 3.1 Based posets \endsubhead

\definition{Definition 3.1.1}
A {\it based poset} is a triple $(\Lambda, \leq, f)$ where 
$(\Lambda, \leq)$ is a poset and
$f$ is a function with the property that for each $\l \in \Lambda$, 
$f(\l)$ is a based ring.  An isomorphism of based posets $$\a :
(\Lambda, \leq, f) \ra (\Pi, \leq', g)$$ is an 
isomorphism of posets
$\a : (\Lambda, \leq) \ra (\Pi, \leq')$ such that for all $\l \in
\Lambda$, $\a$ induces an isomorphism of 
based rings $\a : f(\l) \ra g(\a(\l))$.
\enddefinition

\definition{Definition 3.1.2}
Let $A$ be a reduced tabular algebra with trace and with tabular basis
$\BB$.  Let $(\Lambda, \Gamma, B, M, C, *)$ be a table datum for $A$.

The {\it based poset $P(A, \BB) = P(\Lambda, \Gamma, B, M, C, *) = 
P(\Lambda, \Gamma, B, M, *)$ associated to the table datum $(\Lambda,
\Gamma, B, M, C, *)$}
is the triple $(\Lambda, \leq, f)$ where $(\Lambda, \leq)$ is the
poset in the table datum of $A$ and, for each $\l \in \Lambda$,
$f(\l)$ is the based ring $M_{|M(\l)|, \Gamma(\l), B(\l)}
\cong A_\l^\infty$.
Such a based poset is equipped with an anti-automorphism, $*$, which by
definition leaves the elements of the poset fixed and induces the map
$*$ of Lemma 1.2.4 (ii) on each $f(\l)$.

The category ${\Cal P}(A, \BB)$ has as objects all elements $P(D)$ for
$D \in {\Cal D}(A, \BB)$ (see Definition 2.4.1); the morphisms are
isomorphisms of based posets.
\enddefinition

\definition{Definition 3.1.3}
Let $A$ be a reduced tabular algebra with trace and tabular basis
$\BB$.  Let $D$ be an object of ${\Cal D}(A, \BB)$ and let $P(D)$ be
the corresponding object of ${\Cal P}(A, \BB)$.  A {\it
parametrization} of $D$ consists of 
bijections $$s_\l : M(\l) \ra \{1, 2, \ldots, |M(\l)|\}$$ and a
map $$
\pi : 
\coprod_\l \left( M(\l) \times B(\l) \times M(\l) \right) \ra 
\coprod_{\l \in \Lambda} f(\l)
$$ such that for all $\l \in \Lambda$, $S, T \in M(\l)$ and $b \in B(\l)$,
$\pi(S, b, T) = e_{s_\l(S), s_\l(T)} \otimes b \in f(\l)$.  We will typically
refer to the parametrization $(\pi, \ \coprod_{\l \in \Lambda} s_\l)$ as
``the parametrization $\pi$'' for short.
\enddefinition

\proclaim{Proposition 3.1.4}
Let $A$ be a reduced tabular algebra with trace and tabular basis
$\BB$.  Let $D_1 := (\Lambda, \Gamma, B, M, *)$ and $D_2 := (\Lambda',
\Gamma', B', M', *')$ be objects of ${\Cal D}(A, \BB)$ with
parametrizations $\pi_1$ and $\pi_2$ respectively.
If there is a morphism $\th : D_1 \ra D_2$
in ${\Cal D}(A, \BB)$, then there is a morphism $$
P(\th) : P(D_1) := (\Lambda, \leq, f)
\ra P(D_2) := (\Lambda', \leq', f')
$$ in ${\Cal P}(A, \BB)$, depending on $\pi_1$ and $\pi_2$.  In particular, the
isomorphism type of $P(A, \BB)$ as a based poset is independent of the
choice of table datum.
\endproclaim

\demo{Proof}
Proposition 2.3.3 shows that the sets $\bc_\l$ are independent of the
table datum, so that $$
\th \left( M(\l) \times B(\l) \times M(\l) \right) = M'(\l')
\times B'(\l') \times M'(\l')
$$ for some $\l'$ depending on $\l$.  Since $A$ is reduced, Proposition
2.3.3 also shows that $\th$ is compatible with the partial orders on the
two table data, and thus that $\th$ is a poset isomorphism.

We define the map $P(\th) : \coprod f(\l) \ra \coprod
g(\th(\l))$ (depending on $\pi_1$ and $\pi_2$) by $\zed$-linear
extension of the condition that $P(\th) \circ \pi_1 = \pi_2 \circ \th$.
The map $P(\l)$ respects the partitions induced by $\Lambda$ and
$\Lambda'$.  It is an isomorphism
of based posets by Corollary 2.3.4 and Lemma 2.2.2 (i), because it 
represents the identity map on $A^{\infty}$ with respect to certain bases.
The conclusion follows.
\qed\enddemo

The involution $*$ is respected by the map $P$ in the following sense.

\proclaim{Lemma 3.1.5}
Let $A$ and $\BB$ be as in Proposition 3.1.4, and let $D := (\Lambda,
\Gamma, B, M, *)$ be an object of ${\Cal D}(A, \BB).$  
Then the involution $*$ determines
and is determined by the maps $*$ on the based poset $P(D)$ in ${\Cal
P}(A, \BB)$ given in Definition 3.1.2, and furthermore, the correspondence is
independent of the parametrization chosen.
\endproclaim

\demo{Proof}
The map $*$ induces an obvious permutation of each set $M(\l) \times
B(\l) \times M(\l)$ for each $\l \in \Lambda$.  Choose a
parametrization $\pi$ of $D$ and define the map $P(*) : \coprod
f(\l) \ra \coprod f(\l)$ by the condition $P(*) \circ \pi = \pi \circ *$.
The map $P(*)$ is equal to the based poset anti-automorphism $*$ of Definition 
3.1.2, and the map is independent of the parametrization chosen.  The
converse is easily checked: the map $*$ of $D$ may be reconstructed from the
maps $*$ on the based posets $P(D)$, again independently
of the choice of parametrization.
\qed\enddemo

Because of Lemma 3.1.5, we may identify the map $*$ of $D$ with the
anti-iso\-morph\-ism of the based poset $P(D)$, and we may denote them
both by $*$.

Proposition 3.1.4 hints that $P$ may be a functor, which will turn out
to be the case (see Theorem 3.4.1 (iv)).
The raison d'\^etre of based posets is the following result, which can
be regarded as the main result of the paper and the answer to Question 2.1.4.

\proclaim{Theorem 3.1.6}
Let $A$ be a reduced tabular algebra with trace and with tabular basis
$\BB$.  The categories ${\Cal D}(A, \BB)$ and ${\Cal P}(A, \BB)$ are
equivalent.
\endproclaim

\subhead 3.2 Classifying the isomorphisms of based posets \endsubhead

To understand the morphisms in the category ${\Cal P}(A, \BB)$, we
require the following definition.

\definition{Definition 3.2.1}
Let $(\Lambda, \leq, f)$ and $(\Pi, \leq', g)$ be objects of 
${\Cal P}(A, \BB)$.  Let $p_\Lambda : \Lambda \ra \Pi$ be an isomorphism of 
abstract posets, and for each $\l \in \Lambda$ let $p_\l$ be an 
untwisted isomorphism of based rings (see
Definition 1.3.2) from $f(\l)$ to $g(p_\Lambda(\l))$.
We define the isomorphism of based posets $$\iota_p : 
(\Lambda, \leq, f) \ra
(\Pi, \leq', g)$$ to be the isomorphism inducing the map $p_\Lambda$ on 
$\Lambda$ and the maps $p_\l$ on each $f(\l)$.

Let $(\Lambda, \leq, f)$ be an object of ${\Cal P}(A, \BB)$,
let $\l \in \Lambda$ and set $n = |M(\l)|$.  
Let $\be(b_1, b_2, \ldots, b_n)$, $\s_w$
and $\psi_\e$ be based ring automorphisms of 
$M_{n, \Gamma(\l), B(\l)}$.  We define the based poset automorphism
$\be^\l(b_1, b_2, \ldots, b_n)$ (respectively, $\s^\l_w$, $\psi^\l_\e$) of
$(\Lambda, \leq, f)$ to be the automorphism that induces the identity
map on the underlying poset and on all based rings 
$f(\mu)$ for $\mu \ne \l$, and that induces the automorphism 
$\be(b_1, b_2, \ldots, b_n)$ (respectively, $\s_w$, $\psi_\e$) on $f(\l)$.
\enddefinition

\remark{Remark 3.2.2}
Definition 3.2.1 makes sense by Lemma 2.2.2 (i), which guarantees that the
based rings involved are isomorphic to $M_{n, \Gamma, B}$ for suitable 
$n$, $\Gamma$ and $B$.
\endremark


\proclaim{Proposition 3.2.3}
Maintain the notation of Definition 3.2.1. 
Any morphism $$\a : (\Lambda, \leq, f) \ra (\Pi, \leq', g)$$ in 
${\Cal P}(A, \BB)$ 
can be expressed as a product of isomorphisms of the form $\iota_p$ and
$\be^\l(b_1, b_2, \ldots, b_n)$.
\endproclaim

\demo{Proof}
Corollary 1.3.3 reduces the problem to the case where 
$(\Lambda, \leq, f) = (\Pi, \leq', g)$, once we compose $\a$ with a suitable
isomorphism $\iota_p$.  The result now follows from Proposition 1.3.5,
because the automorphisms $\s_w$ and $\psi_\e$ are untwisted.
\qed\enddemo

\proclaim{Proposition 3.2.4}
Let $A$ be a tabular algebra with trace and tabular basis $\BB$, and
let $$
\a : (\Lambda, \leq, f) \ra (\Lambda', \leq', f')
$$ be a morphism in ${\Cal P}(A, \BB)$.  Then $\a$ intertwines the
based poset anti-automorphisms of its source and target, and $\a$
takes distinguished idempotents to distinguished idempotents.
\endproclaim

\demo{Proof}
By Proposition 3.2.3, it is enough to verify this for $\a = \iota_p$ and
$\a = \be^\l(b_1, b_2, \ldots, b_n)$.  The first case follows easily from
the definitions, and the second case is a consequence of Lemma 1.2.4.
\qed\enddemo

\subhead 3.3 Automorphisms of the tabular basis \endsubhead

It will turn out that morphisms in the category ${\Cal D}(A, \BB)$ all
arise from the following construction.

\definition{Definition 3.3.1}
Let $A$ be a reduced tabular algebra with trace and tabular basis
$\BB$.  Let $D_1 := (\Lambda, \Gamma, B, M, *)$ and 
$D_2 := (\Lambda', \Gamma', B', M', *')$ be objects of ${\Cal D}(A,
\BB)$, with parametrizations $\pi_1$ and $\pi_2$ respectively.

Let $\a : P(D_1) \ra P(D_2)$ be a morphism in ${\Cal P}(A, \BB)$; such
a morphism exists by Proposition 3.1.4.  
This induces a map $$
\a(\pi_1, \pi_2) :
\coprod_{\l \in \Lambda} \left( M(\l) \times B(\l) \times M(\l)
\right) \ra 
\coprod_{\l' \in \Lambda'} \left( M'(\l') \times B'(\l') \times M'(\l')
\right)
$$  given by $\a(\pi_1, \pi_2) = \pi_2^{-1} \a \pi_1$.  
If maps $C$ and $C'$ are chosen such that $$(\Lambda, \Gamma, B,
M, C, *)$$ and $$(\Lambda', \Gamma', B', M', C', *')$$ are table data,
$\a(\pi_1, \pi_2)$ induces a permutation $\a$ of $\BB$ via $$
\a := C' \circ \a(\pi_1, \pi_2) \circ C^{-1}
.$$  This may be extended $\A$-linearly to a map on $A$.
\enddefinition

\remark{Remark 3.3.2}
It must be emphasised that the permutations of the tabular basis in
Definition 3.3.1 are generally {\sl not} algebra automorphisms of $A$.
\endremark

\proclaim{Lemma 3.3.3}
Maintain the notation of Definition 3.3.1, so that the morphism $\a$
in ${\Cal P}(A, \BB)$ induces a permutation of $\BB = \Im(C) = 
\Im(C')$.  Then for each $\l \in \Lambda$,
$\a$ induces a bijection from $M(\l)$ to $M'(\a(\l))$: if $S \in
M(\l)$, $\a(S)$ is defined by the condition
$\a(\pi_1, \pi_2)(S, 1, S) = (\a(S), 1, \a(S)).$  In particular, $\a$ fixes
the distinguished involutions setwise.
\endproclaim

\demo{Proof}
The map $\a$ is compatible with the partitions of $\BB$ by the sets
$\Lambda$ and $\Lambda'$ because it is a morphism of based posets.
The bijections $\pi_1 C^{-1}$ and $\pi_2 (C')^{-1}$ send distinguished
involutions to distinguished idempotents, so by
Proposition 3.2.4 the permutation of $\BB$ induced by $\a$ takes
distinguished involutions to distinguished involutions.
\qed\enddemo

What is remarkable about Definition 3.3.1 is that the
permutations of $\BB$ arising are in fact morphisms in the category
${\Cal D}(A, \BB)$, and that these morphisms can be made to act on the
possible table data.

\proclaim{Proposition 3.3.4}
Let $A$ be a reduced tabular algebra with trace and with tabular basis
$\BB$.  Let $D_1 = (\Lambda, \Gamma, B, M, *)$ and $D_2 = 
(\Lambda', \Gamma', B', M', *')$  be objects of ${\Cal
D}(A, \BB)$ and let $C$ be a map such that $(\Lambda, \Gamma, B, M, C,
*)$ is a table datum.  Fix parametrizations $\pi_1$ and $\pi_2$ for $D$ and 
$D'$ respectively.  Let
$\a$ be a morphism $P(D) \ra P(D')$; this exists by Proposition 3.1.4.
Then $\a(\pi_1, \pi_2)$ is a morphism from $D \ra D'$ in 
${\Cal D}(A, \BB)$ and 
$(\Lambda', \Gamma', B', M', \a(C), *')$ is another reduced table datum for
$A$ (with respect to the same trace) where $\a(C) := 
C \circ \a(\pi_1, \pi_2)^{-1}$.  Furthermore, $*$ and $*'$ are equal as
permutations of $\BB$.
\endproclaim

\demo{Proof}
The last assertion follows from 
Lemma 1.2.4 (iii), Lemma 2.2.2 (i) and Lemma 3.1.5.  For the other assertion,
we check each of the five axioms.

\item{A1.}{It is clear that $\a(C)$ is injective because $\a(\pi_1, \pi_2)$ is
invertible and $C$ is injective.  The other assertions
follow easily from the definition of ${\Cal P}(A, \BB)$.}


\item{A2.}{Proposition 3.2.4 and Lemma 3.1.5 show that 
$\a(\pi_1, \pi_2)^{-1} \circ *' = * \circ \a(\pi_1, \pi_2)^{-1}$.
The map $C$ intertwines the maps $*$ on its domain and range by
axiom (A2) applied to $C$.
We therefore have $$\eqalign{
(C \circ \a(\pi_1, \pi_2)^{-1}(S, b, T))^{*'}
&= (C \circ \a(\pi_1, \pi_2)^{-1}(S, b, T))^* \cr
&= C \circ ((\a(\pi_1, \pi_2)^{-1}(S, b, T))^*) \cr
&= C \circ (\a(\pi_1, \pi_2)^{-1}(S, b, T)^{*'}) \cr
&= C \circ \a(\pi_1, \pi_2)^{-1}(T, \bar{b}, S) \cr
}$$ as required.}

\item{A3.}{We note that $\a$, being an isomorphism of based posets,
respects both the partition of $\BB$ into subsets $\bc_\l$ and 
the partial order on $\Lambda$.  The claims regarding $A(< \l)$ follow
from this.  In particular, $\a$ induces a 
bijection of $\Lambda \ra \Lambda'$, also denoted by $\a$.

We now need to show the existence of a function $r'_a$ with the
required properties with respect to the candidate $(\Lambda, \Gamma,
B, M, \a(C), *)$ for the cell datum.  We write $r_a$ for the
corresponding function associated to the original cell datum.

We need only check the cases $\a = \iota_p$ and 
$\a = \be^\l(b_1, b_2, \ldots, b_n)$ by Proposition 3.2.3.

For $\a = \iota_p$ we have $$
\a(C)(S, \ b, \ T) = C(\s(S), \ \psi(b), \ \s(T))
$$ for all $\l \in \Lambda'$, $S, T \in M'(\l)$ and $b \in B'(\l)$.
Here, $\s : M'(\l) \ra M(\a^{-1}(\l))$ is the map given in Lemma 3.3.3 and 
$\psi : (\Gamma', B') \ra (\Gamma, B)$ is the isomorphism of based rings
determined by $p$.  In this case, we define $$
r'_a(S', S) := \psi^{-1}(r_a(\s(S'), \s(S)))
.$$  Axiom (A3) applied to $C$ now gives $$\eqalign{
a . C_{\s(S), \s(T)}^{\psi(b)} 
&= \sum_{\s(S') \in M(\a^{-1}(\l))} 
C_{\s(S'), \s(T)}^{r_a(\s(S'), \s(S)) \psi(b)} \mod A(< \a^{-1}(\l)) \cr
&= \sum_{\s(S') \in M(\a^{-1}(\l))}
C_{\s(S'), \s(T)}^{\psi(r'_a(S', S) b)} \mod A(< \a^{-1}(\l)). \cr
}$$  This yields $$
a . \a(C)_{S, T}^b = \sum_{S' \in M'(\l)} \a(C)_{S', T}^{r'_a(S', S)b}
\mod A(<' \l)
$$ and shows that the axiom holds with respect to the prospective table datum 
for $D'$ and $r'_a$ in place of $r_a$.

For the other case, we take $\a = \be^\l(b_1, b_2, \ldots, b_n)$.  In
this case, $D = D'$ and so $\pi_1 = \pi_2 = \pi$. 
For each $S \in
M(\l)$, let us write $b_S$ for $b_{s_\l(S)}$, where $s_\l$ is 
associated to the parametrization $\pi$ in the usual way.
In this case, $$
\a(C)(S, \ b, \ T) := C(S, \ \overline{b_S} b b_T, \ T)
.$$ (Recall that $\overline{b_S} b b_T \in B(\l)$ by Lemma 1.2.2 
because $b_S$ and $b_T$ are grouplike.)  In this case, we set
$r'_a(S', S) := b_{S'} r_a(S', S) \overline{b_S}$.
Axiom (A3) applied to $C$ now gives $$\eqalign{
a . C_{S, T}^{\overline{b_S} b b_T}
&= \sum_{S' \in M(\l)} 
C_{S', T}^{r_a(S', S) \overline{b_S} b b_T} \mod A(< (\l)) \cr
&= \sum_{S' \in M(\l)} 
C_{S', T}^{\overline{b_{S'}} r'_a(S', S) b b_T} \mod A(< (\l)). \cr
}$$ This yields $$
a . \a(C)_{S, T}^b = \sum_{S' \in M'(\l)} \a(C)_{S', T}^{r'_a(S', S)b}
\mod A(< \l)
$$ as required.}

\item{A4.}{As in the verification of axiom (A3), we need only check the
cases $\a = \iota_p$ and $\a = \be^\l(b_1, b_2, \ldots, b_n)$.  The case of
$\a = \iota_p$ is a routine verification.  
For the other case, the condition
for the degree bound to be achieved follows from the observation
that $$\k(\overline{b_S} b'' b_V, (\overline{b_S} b b_T)
(\overline{b_T} b' b_V)) = \k(b'', bb')
.$$  The claim about the case $b = b' = b'' = 1$ follows from Lemma 
2.2.2 (i) and the fact
that the product $$(e_{s_\l(S), s_\l(T)} \otimes 
\overline{b_S} b_T)(e_{s_\l(T), s_\l(V)} \otimes \overline{b_T} b_V)$$
contains $e_{s_\l(S), s_\l(V)} \otimes \overline{b_S} b_V$ with
coefficient $1$.}

\item{A5.}{The map $C \circ \a(\pi_1, \pi_2) \circ C^{-1}$ 
sends distinguished involutions to distinguished involutions by Lemma
3.3.3, and it commutes with the map $* = *'$ on $\BB$ by the argument
establishing axiom (A2) above.  Axiom (A5) follows easily from these
observations.}
\qed\enddemo

\subhead 3.4 Main results \endsubhead

We are now in a position to examine the category ${\Cal D}(A, \BB)$.

\proclaim{Theorem 3.4.1}
Let $A$ be a reduced tabular algebra with trace and tabular basis
$\BB$.  Let $D_1 := (\Lambda, \Gamma, B, M, *)$ and 
$D_2 := (\Lambda', \Gamma', B', M', *')$ be objects of ${\Cal D}(A,
\BB)$, with parametrizations $\pi_1$ and $\pi_2$ respectively.

\item{\rm (i)}{Let $\phi : P(D_1) \ra P(D_2)$ be a morphism in ${\Cal
P}(A, \BB)$ and choose $C$ so that $(\Lambda, \Gamma, B, M, C, *)$
is a table datum.  Then $\phi(\pi_1, \pi_2) : D_1 \ra D_2$ is a
morphism in ${\Cal D}(A, \BB)$ and there exists $C'$ such that
$(\Lambda', \Gamma', B', M', C', *')$ is a cell datum and such that
$\phi(\pi_1, \pi_2) = (C')^{-1} \circ C$.}
\item{\rm (ii)}{Every morphism $\phi : D_1 \ra D_2$ in ${\Cal D}(A,
 \BB)$ is of the form $\a(\pi_1, \pi_2)$ for some morphism $\a : P(D_1)
\ra P(D_2)$ in ${\Cal P}(A, \BB)$.}
\item{\rm (iii)}{The category ${\Cal D}(A, \BB)$ is well defined.}
\item{\rm (iv)}{The map $P$ is a functor from ${\Cal D}(A, \BB)$ to
${\Cal P}(A, \BB)$ (assuming each object of ${\Cal D}(A, \BB)$ is
assigned a parametrization).}
\endproclaim

\demo{Proof}
Part (i) is immediate from Proposition 3.3.4.

To prove (ii), we first fix such a $\phi$.  Proposition 3.1.4 produces
a morphism $$P(\phi) : P(D_1) \ra P(D_2)$$ in ${\Cal P}(A, \BB)$ 
depending on $\pi_1$ and $\pi_2$.  Definition 3.3.1 then shows that 
$P(\phi)(\pi_1, \pi_2) = \phi$.

For (iii), let $\phi : D_1 \ra D_2$ be a morphism in ${\Cal D}(A,
\BB)$ and choose $C$ as in the statement of (i).  By (ii), $\phi$ is
of the form $\a(\pi_1, \pi_2)$ for some morphism $\a$ in ${\Cal P}(A,
\BB)$.  Applying (i) to $\a$, we see that $\phi = (C')^{-1} \circ C$,
where $C'$ is such that $(\Lambda', \Gamma', B', M', C', *')$ is a
cell datum.  This shows that the morphisms defined in Definition 2.4.1
do not depend on the choice of $C$.  It also shows that composition of 
morphisms is 
well-defined, because if $\phi = (C_1)^{-1} C_2$ and $\phi' =
(C_3)^{-1} C_4$ are morphisms in ${\Cal D}(A, \BB)$, we may arrange
for $C_2 = C_3$ so that $\phi \circ \phi'$ is a morphism.  This proves
(iii), and makes (iv) into an easy exercise.
\qed\enddemo

We can now prove the main result.

\demo{Proof of Theorem 3.1.6}
It is enough to prove that $P$ is an equivalence of categories.

Consider the full subcategory ${\Cal D}'(A, \BB)$ of ${\Cal D}(A,
\BB)$ whose objects are those $(\Lambda, \Gamma, B, M, *)$ for which 
each set $M(\l)$ consists
of the first $|M(\l)|$ natural numbers.  This object has a
parametrization in which all the maps $s_\l$ are the identity.
Furthermore, each object $X$ of ${\Cal P}(A, \BB)$ arises as $P(D')$ for a
unique $D' \in {\Cal D'}(A, \BB)$.  In this case, we define $Q(X) := D'$.  If
$\a : X \ra Y$ is a morphism in ${\Cal P}(A, \BB)$, we define the morphism
$Q(\a) : Q(X) \ra Q(Y)$ by $Q(\a) = \a(\pi, \pi)$, where $\pi$ is the
identity parametrization.  Theorem 3.4.1 shows that $Q$ is a functor
and that $P Q$ is the identity functor on ${\Cal P}(A, \BB)$.

Suppose all objects of ${\Cal D}(A, \BB)$ have been assigned parametrizations.
Let $$D_1 := (\Lambda, \Gamma, B, M, *)$$ be an object of 
${\Cal D}(A, \BB)$.  Let us write $$
(\Lambda, \Gamma, B, M', *) := Q P(D_1)
.$$  The parametrization $\pi_1$ of $D_1$ induces a morphism (\idest
an isomorphism) $\eta_1$ in ${\Cal
D}(A, \BB)$ from $D_1$ to $Q P(D_1)$: this is the map that sends 
$(S, b, T)$ to $(s_\l(S), b, s_\l(T))$ for each $\l \in \Lambda$, 
$S, T \in M(\l)$, 
$b \in B(\l)$ and $s_\l(S), s_\l(T) \in M'(\l)$.  If $D_2$ is another
object with $\phi : D_1 \ra D_2$ a morphism in ${\Cal D}(A, \BB)$, we
then see that the maps $\eta$ give natural isomorphisms between the
identity functor on ${\Cal D}(A, \BB)$ and the functor $Q P$.
Theorem 3.1.6 follows.
\qed\enddemo

\head 4. Algebra automorphisms of tabular algebras \endhead

Recall that in Remark 3.3.2,  we stated that permutations of the
tabular basis arising from morphisms in ${\Cal D}(A, \BB)$
do not always induce isomorphisms at the level of
tabular algebras.  However, the permutations of the tabular
basis that do give algebra automorphisms may be understood using our
results, and may be classified in natural examples.  We do this for
the Brauer algebra in \S4.2.  We do not claim that these results are original;
they are provided merely to illustrate the results of this paper.

\proclaim{Proposition 4.1.1}
Let $A$ be a reduced tabular algebra with trace and with tabular basis
$\BB$.  Let $(\Lambda, \Gamma, B, M, C, *)$ be a table datum for $A$.
Let $\phi$ be an $\A$-algebra automorphism of $A$
preserving $\BB$ setwise.  Then $\phi$ is of the form $C \circ \th
\circ C^{-1}$ for some $\th \in \aut_{{\Cal D}(A, \BB)}(D)$, where $D
:= (\Lambda, \Gamma, B, M, *)$.
\endproclaim

\demo{Proof}
Since $\phi$ is an algebra automorphism preserving the tabular basis,
it must (by Proposition 2.3.3) 
permute the collection of sets $\bc_\l$, so $\phi$
induces a bijection $\a : \Lambda \ra \Lambda$.   This bijection is an
isomorphism of posets because the tabular algebra is reduced and thus
the partial order is determined by the algebra structure via axiom (A3).
It is clear from the definition of $\afn$ that $\afn(\a(X)) = \afn(X)$ 
for all $X \in \Im(C)$. 
Corollary 2.3.4 shows that $\phi$ induces an isomorphism of based rings $$
\a : \coprod_\l A_\l^\infty \ra \coprod_\l A_\l^\infty
.$$  If we choose a parametrization $\pi$ for $D$ and let $J$ be the
map sending $X \in \BB$ to $t_X \in A^\infty$, we find that $$
\pi C^{-1} J^{-1} : 
\coprod A_{\l}^\infty \ra \coprod M_{|M(\l)|, \Gamma(\l), B(\l)}
$$ extends to give an isomorphism between $P(D)$ and the based poset
under consideration.  The based poset isomorphism $\a$ thus gives a
based poset isomorphism $\pi C^{-1} J^{-1} \a J C \pi^{-1} \in
\aut_{{\Cal D}(A, \BB)}(P(D))$.  Since $P$ is an equivalence of
categories by the proof of Theorem 3.1.6, there is a morphism $\th \in 
\aut_{{\Cal D}(A, \BB)}(D)$ such that $\th = C^{-1} J^{-1} \a J C$.
The claim follows from the fact that $\phi = J^{-1} \a J$.
\qed\enddemo

We can deduce the following result, which holds for any table datum.

\proclaim{Corollary 4.1.2}
Let $A$ be a reduced tabular algebra with trace and with tabular basis
$\BB$.  Let $\phi$ be an $\A$-algebra automorphism of $A$
preserving $\BB$ setwise.  Then $\phi$ is a $*$-automorphism 
(\idest $\phi \circ *$ = $* \circ \phi$) and
$\phi$ takes distinguished involutions to distinguished involutions.
\endproclaim

\demo{Proof}
This follows from Proposition 4.1.1, Proposition 3.2.4 and the 
equivalence of categories given by Theorem 3.1.6.
\qed\enddemo

\subhead 4.2 The Brauer algebra \endsubhead

We now recall how Brauer's centralizer algebra (which we call ``the
Brauer algebra'' for short) may be given the structure of a
tabular algebra with trace and show how the results and techniques of
this paper apply to it.  We calculate the group of algebra automorphisms of
the algebra that preserve the basis; many other natural examples of
tabular algebras can be analysed similarly.

Some useful references on the Brauer algebra are \cite{{\bf 2}},
\cite{{\bf 3}, \S4} and \cite{{\bf 12}}.

\definition{Definition 4.2.1}
The Brauer algebra $B_n$ ($n \geq 2$)
is defined to be the unital associative
$\A$-algebra with basis consisting of certain diagrams.  A basis
diagram, $D$, consists of two rows of $n$ points, labelled $\{1,
\ldots, n\}$, with each point joined to precisely one point distinct
from itself.  (See Figure 1.)


\topcaption{Figure 1} A Brauer algebra basis element for $n = 6$ \endcaption
\centerline{
\hbox to 1.958in{
\vbox to 0.750in{\vfill
        \includegraphics{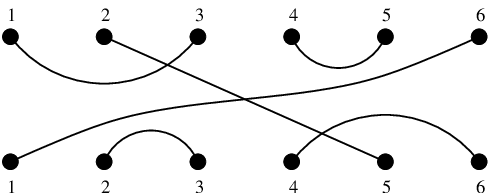}
}
\hfill}
}

Two diagrams $D_1, D_2$ may be composed to obtain $D_1 \circ D_2$ by
placing $D_1$ above $D_2$ and joining corresponding points.  This
produces a number $n(D_1, D_2)$ of interior loops, which are deleted.
The product $D_1 D_2$ in the Brauer algebra is defined by $$
D_1 D_2 := (v + v^{-1})^{n(D_1, D_2)} D_1 \circ D_2
.$$
\enddefinition

As in \cite{{\bf 3}, \S4}, we may describe the basis diagrams in terms of
certain triples.  

\definition{Definition 4.2.2}
Fix a diagram $D$. 
The integer $t(D)$ is defined to be the
number of {\it through strings}, \idest strings joining points in
different rows.  The involutions $S_1(D)$, $S_2(D)$ in the symmetric
group $\Sy(n)$ are defined such that $S_i(D)$ interchanges the ends of
the joins between points in row $i$, with $i \in \{1, 2\}$.
Corresponding to these we have subsets $\text{\rm Fix}(S_i(D))$ of $\{1,
\ldots, n\}$, which are the fixed points of the involutions $S_i(D)$.
Finally, we have a permutation $w(D)$ in $\Sy(t)$,
where $t = t(D)$; this is the permutation of $\text{\rm Fix}(S_1(D))$
determined by taking the end points of the through strings (regarded
as joining from row $2$ to row $1$) in the order determined by taking
their starting points in row $2$ in increasing order.  (We consider
$\Sy(0)$ to be the trivial group, in which case $w$ is the identity.)
The diagram $D$ is then determined by the triple $[S_1(D), S_2(D), w(D)]$.
\enddefinition

We now recall a table datum for the Brauer algebra from \cite{{\bf 4},
Example 2.1.2}.

\definition{Definition 4.2.3}
Let $B_n$ be the Brauer algebra (over $\A$) on $n$ strings.  
The algebra has a table datum $(\Lambda, \Gamma, B, M, C, *)$ as follows.

Take $\Lambda$ to be the set of integers $i$ between $0$ and $n$ such that
$n - i$ is even, ordered in the natural way.  
If $\l = 0$, take $(\Gamma(\l), B(\l))$ to be the
trivial one-dimensional table algebra; otherwise, take $\Gamma(\l)$ to
be the group ring $\zed \Sy(\l)$ with basis $B(\l) = \Sy(\l)$ and 
involution $\overline{w}
= w^{-1}$.  Take $M(\l)$ to be the set of involutions on $n$ letters 
with $\l$ fixed points.  Take $C(S_1, w, S_2) = [S_1,
S_2, w]$; $\Im(C)$ contains the identity element.  
The anti-automorphism $*$ sends $[S_1, S_2, w]$ to $[S_2, S_1, w^{-1}]$.
\enddefinition

\remark{Remark 4.2.4}
There exists a tabular trace $\t$ for $B_n$ (see \cite{{\bf 4}, Remark
2.1.3}).  One way to construct such a trace is to define $\t(D)$ as
follows on basis diagrams.  Count the number, $k(D)$, of loops formed when
each point $i$ in row $1$ is joined to point $i$ in row $2$ by a new
string.  Then the
linear map $\t : B_n \ra \A$ such that $\t(D) := v^{-n} (v +
v^{-1})^{k(D)}$ can be shown to be a trace with these properties.  (We
leave this as an exercise.)
\endremark

Although $\Lambda$ is totally ordered in this case, it can easily be
checked that the table datum given is reduced.  

The large supply of grouplike elements in the sets $B(\l)$ ensures that the 
corresponding based poset has plenty of symmetries, so by Theorem
3.4.1 (i), there are many choices for the map $C$ that give 
the same tabular basis of
diagrams.  As a consequence, there is nothing special about the
diagrams parametrized by elements $C_{S, T}^1$; the set of diagrams
that can be so expressed depends very much on the table datum.  Another
manifestation of this ambiguity is
the fact that the definition of $w(D)$ depends on the choice of two
orderings.  However, the set of all diagrams of the form 
$[S_1(D), S_1(D), 1]$ does not depend on the choice of table datum by 
Lemma 3.3.3, because these are the distinguished involutions.

We wish to calculate the group of permutations of the diagram basis
that preserve the algebra structure.  Examples of such permutations
are those which arise from relabelling the points $\{1, \ldots, n\}$
in rows $1$ and $2$ of each diagram by a fixed permutation in $\Sy(n)$.
(Another way to think about this is to conjugate each diagram by a
fixed diagram with $n$ through strings.)  We will show that all basis
preserving algebra automorphisms of $B_n$ are of this form; in particular,
the outer automorphisms of the group $\Sy(6)$ do not extend to
automorphisms of $B_6$.
The elements $e_{a, b}$ and $g_{a, b}$ of the next definition will
play a key role in the proof.


\definition{Definition 4.2.5}
For $1 \leq a < b \leq n$, we define the basis elements $e_{a, b}$ and
$g_{a, b}$ of $B_n$ as follows.  

For the element $e_{a, b}$, point $j$ in row $1$ is
joined to point $j$ in row $2$, unless $j \in \{a, b\}$.  Points $a$
and $b$ in row $i$ (for $i \in \{1, 2\}$) are joined to each other.

For the element $g_{a, b}$, point $j$ in row $1$ is
joined to point $j$ in row $2$, unless $j \in \{a, b\}$.  Point $a$
in row $i$ is joined to point $b$ in row $3 - i$
(for $i \in \{1, 2\}$).
\enddefinition

It is clear that the elements $g_{a, b}$ generate a subalgebra of $B_n$
isomorphic to $\zed(\Sy(n))$, where $g_{a, b}$ corresponds to the
transposition $(a, b)$.  More importantly, we have the following
well-known fact.

\proclaim{Proposition 4.2.6}
The algebra $B_n$ is generated as a unital $\A$-algebra by the
set $\{e_{k, k+1} : 1 \leq k < n\} \cup \{g_{k, k+1} : 1 \leq k < n\}$.
\endproclaim

\demo{Proof}
See \cite{{\bf 12}, Proposition 2.1 (a)}.
\qed\enddemo

The elements $e_{a, b}$ and $g_{a, b}$ may be identified by the
following properties which are independent of the table datum chosen.
By Proposition 3.1.4, it makes sense to refer to the maximal element of
the poset $\Lambda$ as $\l_0$, and to the second maximal element as
$\l_1$.  (Recall that $\Lambda$ is totally ordered, and that $n \geq
2$ so that $|\Lambda| \geq 2$.)

\proclaim{Lemma 4.2.7}
The elements $e_{a, b}$ are precisely the distinguished involutions 
in the set $\bc_{\l_1}$.  The elements $g_{a, b}$ can be characterized
as the only nonidentity elements in $\bc_{\l_0}$ such that there
exists a distinguished involution $e \in \bc_{\l_1}$ with $g_{a, b} e
= e$.  (If this happens, we have $e = e_{a, b}$.)  
These characterizations depend only on the basis, and not on
the table datum.
\endproclaim

\demo{Proof}
Axiom (A5) shows that the set of basis elements that are
distinguished involutions is independent of the table datum.  The
other assertions follow easily from Definition 4.2.3.
\qed\enddemo

We present the following result to illustrate our techniques and to
confirm Remark 3.3.2.

\proclaim{Proposition 4.2.8}
Let $\a : B_n \ra B_n$ be an automorphism of $\A$-algebras
preserving the diagram basis elements.  Then there exists a diagram
$X$ with $n$ through strings such that $\a(z) = X^{-1} z X$ for all $z
\in B_n$.
\endproclaim

\demo{Proof}
By Proposition 4.1.1, $\a = C \circ \th
\circ C^{-1}$ for some $\th \in \aut_{{\Cal D}(A, \BB)}(D)$, where $D =
(\Lambda, \Gamma, B, M, *)$.  By the main results, there is a corresponding
automorphism $P(\th)$ of the based poset $P(D)$ in ${\Cal P}(A, \BB)$.
Since $\Lambda$ is totally ordered, it has no non-trivial
automorphisms as an abstract poset, so $P(\th)$ fixes each poset
element and $\a$ must fix 
$\bc_{\l_0}$ and $\bc_{\l_1}$ setwise.  By Corollary 4.1.2 $\a$ permutes
the distinguished involutions in the set $\bc_{\l_1}$; in other words,
for each $a$ and $b$ with $1 \leq a < b \leq n$ there exist $1 \leq c
< d \leq n$ with $\a(e_{a, b}) = e_{c, d}$.  By Lemma 4.2.7, we must
have $\a(g_{a, b}) = g_{c, d}$.  This determines $\a$ by Proposition
4.2.6.  It remains to show that $\a$ is of the required form.

The map $\a$ induces an isomorphism of the quotient algebra $A/A(<
\l_0)$ because it fixes each $\bc_\l$ setwise.  This algebra is
naturally isomorphic to $\zed \Sy(n)$, and $\a$ induces an automorphism
of $\Sy(n)$ that preserves cycle type.  It follows that the action of
$\a$ on $A/A(< \l_0)$ is effected by conjugation by an element $g \in
\bc_0$ (\idest $\a(z) = g^{-1} z g$).  The automorphism $$
\phi : z \mapsto g \a(z) g^{-1}
$$ of $B_n$ 
preserves the diagram basis and fixes all the elements $g_{a, b}$.  
By Proposition 4.2.6 and Lemma 4.2.7, $\phi$ is the identity, which
shows that $\a$ is conjugation by $g$ and completes the proof.
\qed\enddemo

\head Acknowledgement \endhead

The author is grateful to the referee for helpful suggestions.

\leftheadtext{}
\rightheadtext{}

\end